\documentclass[letterpaper, 10pt, conference]{ieeeconf}
\IEEEoverridecommandlockouts 
\usepackage{amsmath,amssymb,color}
\usepackage{graphicx,subfigure,warmread}
\usepackage[all,import]{xy}

\newcommand{\norm}[1]{\ensuremath{\left\| #1 \right\|}}

\newcommand{\bracket}[1]{\ensuremath{\left[ #1 \right]}}

\newcommand{\parenth}[1]{\ensuremath{\left( #1 \right)}}

\newcommand{\refeqn}[1]{(\ref{eqn:#1})}
\newcommand{\reffig}[1]{Fig. \ref{fig:#1}}
\newcommand{\tr}[1]{\mbox{tr}\ensuremath{\negthickspace\bracket{#1}}}
\newcommand{\deriv}[2]{\ensuremath{\frac{\partial #1}{\partial #2}}}
\newcommand{\G}{\ensuremath{\mathsf{G}}}
\newcommand{\SO}{\ensuremath{\mathsf{SO(3)}}}
\newcommand{\T}{\ensuremath{\mathsf{T}}}
\renewcommand{\L}{\ensuremath{\mathsf{L}}}
\newcommand{\so}{\ensuremath{\mathfrak{so}(3)}}

\renewcommand{\Re}{\ensuremath{\mathbb{R}}}

\newcommand{\D}{\ensuremath{\mathbf{D}}}

\newcommand{\Ad}{\ensuremath{\mathrm{Ad}}}

\newcommand{\g}{\ensuremath{\mathfrak{g}}}

\renewcommand{\sb}{\ensuremath{\overline{s}}}
\newcommand{\mub}{\ensuremath{\overline{\mu}}}

\title{\LARGE \bf
Dynamics of a 3D Elastic String Pendulum}

\author{Taeyoung Lee, Melvin Leok\authorrefmark{1}, and N. Harris McClamroch\authorrefmark{2}%
\thanks{Taeyoung Lee, Mechanical and Aerospace Engineering, Florida Institute of Technology, Melbourne, FL 39201 {\tt taeyoung@fit.edu}}%
\thanks{Melvin Leok, Mathematics, Purdue University, West Lafayette, IN 47907 {\tt mleok@math.purdue.edu}}%
\thanks{N. Harris McClamroch, Aerospace Engineering, University of Michigan, Ann Arbor, MI 48109 {\tt
nhm@umich.edu}}%
\thanks{\textsuperscript{\footnotesize\ensuremath{*}}This research has been supported in part by NSF under grants DMS-0714223, DMS-0726263, DMS-0747659.}
\thanks{\textsuperscript{\footnotesize\ensuremath{\dagger}}This research has been supported in part by NSF under grant CMS-0555797.}
}

\begin{document}
\allowdisplaybreaks
\maketitle \thispagestyle{empty} \pagestyle{empty}

\begin{abstract}
This paper presents an analytical model and a geometric numerical integrator for a rigid body connected to an elastic string, acting under a gravitational potential. Since the point where the string is attached to the rigid body is displaced from the center of mass of the rigid body, there exist nonlinear coupling effects between the string deformation and the rigid body dynamics. A geometric numerical integrator, refereed to as a Lie group variational integrator, is developed to numerically preserve the Hamiltonian structure of the presented model and its Lie group configuration manifold. These properties are illustrated by a numerical simulation.
\end{abstract}

\section{Introduction}

The dynamics of a body connected to a string appear in several engineering problems such as cable cranes, towed underwater vehicles, and tethered spacecraft. It has been shown that gravitational forces acting along a string can alter the tension of the string, while significantly disturbing the dynamics of a body connected to the string~\cite{McLRocPI92}. Therefore, it is important to model the string dynamics accurately as well as the dynamics of the body even if the tension of the string is low. 

Several dynamic and numerical models have been developed. Lumped mass models, where the string is spatially discretized into connected point masses, were developed in~\cite{ChaOE84,DriLueAOR00,WalPolMC60}. Finite difference methods in both the spatial domain and the time domain were applied in~\cite{KohZhaJEM99,KuhSeiAMC95}. Finite element discretizations of the weak form of the equations of motion were applied in~\cite{KuhSeiAMC95,BucDriTA04}. String deployment models were developed in~\cite{BanKanJAS82,SchSteISDSTA00}.

The goal of this paper is to develop an analytical model and a numerical simulation tool for a rigid body connected to a string acting under a gravitational potential. This dynamic model is referred to as a 3D elastic string pendulum: it is a generalization of a 3D pendulum model introduced in~\cite{SheSanPICDC04} to include the effects of string deformations; it is an extension of a string pendulum model with a point mass bob~\cite{KuhSeiAMC95}. 

We assume that the point where the string is attached to the rigid body is displaced from the center of mass of the rigid body so that there exist nonlinear coupling effects between the string deformation dynamics and the rigid body dynamics. This provides a more realistic and accurate dynamic model. We show that the governing equations of motion can be developed according to Hamilton's variational principle. 

The second part of this paper deals with a geometric numerical integrator for the 3D elastic string pendulum. Geometric numerical integration is concerned with developing numerical integrators that preserve geometric features of a system, such as invariants, symmetry, and reversibility~\cite{HaiLub00}. For a numerical simulation of Hamiltonian systems evolving on a Lie group to exhibit good long-time energy behavior, it is critical to preserve both the symplectic property of Hamiltonian flows and the Lie group structure~\cite{LeeLeoCMDA07}. A geometric numerical integrator, referred to as a Lie group variational integrator, has been developed for a Hamiltonian system on an arbitrary Lie group in~\cite{Lee08}.

A 3D elastic string pendulum is a Hamiltonian system, and its configuration manifold is expressed as the product of the special Euclidean group $\SO$ and the space of connected curve segments on $\Re^3$. This paper develops a Lie group variational integrator for a 3D elastic string pendulum based on the results presented in~\cite{Lee08}. The proposed geometric numerical integrator preserves symplecticity and momentum maps, and exhibits desirable energy conservation properties. It also respects the Lie group structure of the configuration manifold, and avoids the singularities and computational complexities associated with the use of local coordinates.

In summary, this paper develops an analytical model and a geometric numerical integrator for a 3D elastic string pendulum. These provide a mathematical model and a reliable numerical simulation tool that characterizes the nonlinear coupling between the string dynamics and the rigid body dynamics accurately. This can be naturally extended to controlled dynamics, and serve as the basis for optimal control algorithms as in~\cite{LeLeMc2008a}.

This paper is organized as follows. A 3D elastic string pendulum  is described in Section~\ref{sec:3DESP}. An analytical model and a Lie group variational integrator are developed in Section~\ref{sec:AM} and in Section~\ref{sec:LGVI}, respectively, followed by a numerical example in Section~\ref{sec:NE}.

\section{3D Elastic String Pendulum}\label{sec:3DESP}

Consider a rigid body that is attached to an elastic string. The other end of the string is fixed to a pivot point. We assume that the rigid body can freely translate and rotate in a three dimensional space, and the string is extensible and flexible. The bending stiffness of the string is not considered as the diameter of the string is assumed to be negligible compared to its length. The point where the string is attached to the rigid body is displaced from the center of mass of the rigid body so that the dynamics of the rigid body is coupled to the string deformations and displacements. 

\begin{figure}\footnotesize\selectfont
$$\begin{xy}
\xyWARMprocessEPS{3DSP}{eps}%%
\xyMarkedImport{}%%
\xyMarkedMathPoints{1-11}
\end{xy}$$
\caption{3D Elastic String Pendulum}\label{fig:3DSP}
\end{figure}
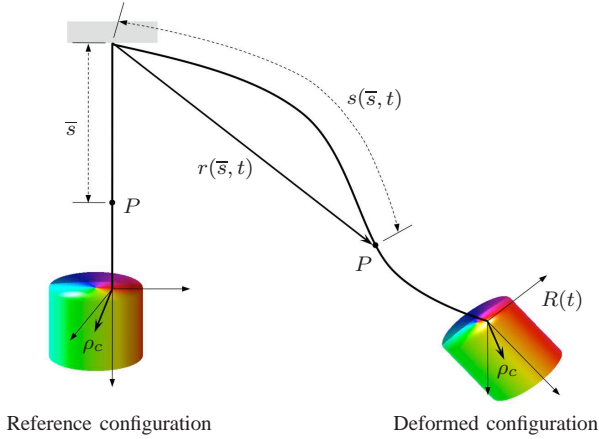

This model is a generalization of the 3D pendulum and the string pendulum introduced in \cite{SheSanPICDC04} and \cite{KuhSeiAMC95}, respectively, and it is referred to as a 3D elastic string pendulum. This is illustrated in \reffig{3DSP}.

\renewcommand{\xyWARMinclude}[1]{\includegraphics[width=0.85\columnwidth]{#1}}

We choose a global reference frame and a body-fixed frame. The origin of the body-fixed frame is located at the end of the string where the string is attached to the rigid body. Since the string is extensible, we need to distinguish between the arc length for the stretched deformed configuration and the arc length for the unstretched reference configuration. Define

\vspace*{0.2cm}
{\allowdisplaybreaks
\centerline{
\begin{tabular}{lp{6.3cm}}
$l\in\Re$ & Total length of the unstretched string\\
$\sb \in[0,l]$ & Length of the string from the pivot to a material point $P$ for the unstretched reference configuration\\
$s(\sb,t) \in\Re^+$ & Length of the string from the pivot to a material point $P$ for the stretched deformed configuration\\
$r(\sb,t)\in\Re^3$ & Vector from the pivot to a material point $P$ in the global reference frame\\
$R\in\SO$ & Rotation matrix from the body-fixed frame to the reference frame\\
$\Omega\in\Re^3$ & Angular velocity of the rigid body represented in the body-fixed frame\\
$\rho_c\in\Re^3$ & Vector from the origin of the body fixed frame to the mass center of the rigid body represented in the body fixed frame\\
$\overline{\mu}\in\Re^+$ & Mass density of the string per unit unstretched length \\
$M\in\Re^+$ & Mass of the rigid body\\
$J\in\Re^{3\times3}$ & Inertia matrix of the rigid body represented in the body fixed frame
\end{tabular}}
\vspace*{0.2cm}

A configuration of this system can be described by the locations of all the material points of the string, $r(\sb,t)$ for $\sb\in[0,l]$, and the attitude of the rigid body $R(t)$ with respect to the reference frame. So, the configuration manifold is $\G=C^\infty([0,l],\Re^3)\times\SO$, where $C^\infty([0,l],\Re^3)$ denotes the space of smooth connected curve segments on $\Re^3$ and $\SO=\{R\in\Re^{3\times 3}\,|\, R^TR=I, \det[R]=1\}$.

\section{Continuous-time Analytical Model}\label{sec:AM}

In this section, we develop continuous-time equations of motion for a 3D elastic string pendulum. The equations for a string pendulum connected to a point mass has been developed in \cite{KuhSeiAMC95}. Here we focus on generalizing them for a rigid body. The attitude kinematics equation of the rigid body is given by
\begin{align}
    \dot R = R\hat\Omega,\label{eqn:Rdot}
\end{align}
where the \textit{hat map} $\hat\cdot:\Re^3\rightarrow\so$ is defined by the condition that $\hat x y=x\times y$ for any $x,y\in\Re^3$. For notational simplicity, we do not express the time dependency of variables explicitly, i.e. $r(\sb)=r(\sb,t)$.

\subsection{Lagrangian}

\paragraph*{Kinetic energy}
The total kinetic energy is composed of the kinetic energy of the string $T_{str}$ and the kinetic energy of the rigid body $T_{rb}$. Let $\dot r(\sb,t)$ be the partial derivative of $r(\sb,t)$ with respect to $t$. This represents the velocity of a material point on the string. Then, the kinetic energy of the string is given by
\begin{align}
T_{str} = \int_{0}^l \frac{1}{2}\mub \norm{\dot r(\sb)}^2 \, d\sb.\label{eqn:Tstr}
\end{align}
Let $\rho\in\Re^3$ be the vector from the mass center of the rigid body to a mass element of the rigid body represented in the body fixed frame. The location of the mass element is given by $r(l)+R(\rho_c+\rho)$ in the global reference frame. Therefore, the kinetic energy of the rigid body can be written as
\begin{align}
T_{rb} & = \int_{\mathcal{B}} \frac{1}{2}\|\dot r(l) + R\hat\Omega(\rho_c+\rho)\|^2\,dM(\rho)\nonumber\\
& = \frac{1}{2}M \dot r(l)\cdot \dot r(l) +\frac{1}{2}\Omega\cdot J\Omega +M\dot r(l)\cdot R\hat\Omega\rho_c,
\end{align} 
where $\mathcal{B}$ denotes the region enclosed by the rigid body surface, and we use the following properties: $\int_{\mathcal{B}}\rho\,dM=0$; $\hat x y = -\hat y x$; $J=-\int_{\mathcal{B}}((\rho+\rho_c)^{\wedge})^2\,dM$.

\paragraph*{Potential Energy}
The strain of the string at a material point located at $r(\sb)$ is given by
\begin{align*}
\epsilon = \lim_{\Delta \sb\rightarrow 0} \frac{\Delta s(\sb)-\Delta \sb}{\Delta \sb}
= s'(\sb) -1,
\end{align*}
where $(\;)'$ denote the partial derivative with respect to $\sb$. The tangent vector at the material point is given by 
\begin{align*}
e_t= \deriv{r(\sb)}{s} = \deriv{r(\sb)}{\sb}\deriv{\sb}{s(\sb)} = \frac{r'(\sb)}{s'(\sb)}.
\end{align*}
Since this tangent vector has the unit length, we have $s'(\sb)=\norm{r'(\sb)}$. Therefore, the strain is given by $\epsilon = \norm{r'(\sb)}-1$. The potential energy of the string is composed of the elastic potential and the gravitational potential:
\begin{align}
V_{str} = \int_{0}^l \frac{1}{2}EA (\norm{r'(\sb)}-1)^2 -\mub g r(\sb)\cdot e_3 \,d\sb,
\end{align}
where $E$ and $A$ denote the Young's modulus and the sectional area of the string, respectively, and the unit vector $e_3$ represents the gravity direction.

Since the location of the center of mass of the rigid body is $r(l)+R\rho_c$ in the global reference frame, the gravitational potential energy of the rigid body is
\begin{align}
V_{rb} = -Mg (r(l)+ R\rho_c) \cdot e_3.\label{eqn:Vrb}
\end{align}

From \refeqn{Tstr}-\refeqn{Vrb}, the Lagrangian of the 3D elastic string pendulum is given by
\begin{align}
L = T_{str}- V_{str} + T_{rb} - V_{rb}.\label{eqn:L}
\end{align}

\subsection{Euler-Lagrange Equations}
Let the action integral be $\mathfrak{G} = \int_{t_0}^{t_f} L\,dt$. It is composed of two parts, $\mathfrak{G}_{str}$ and $\mathfrak{G}_{rb}$, contributed by the string and by the rigid body, respectively. According to the Hamilton's principle, the variation of the action integral is equal to zero for fixed boundary conditions, which yields the Euler-Lagrange equations of the 3D elastic string pendulum.

By repeatedly applying integration by parts, the variation of $\mathfrak{G}_{str}$ can be written as
\begin{align}
\delta&\mathfrak{G}_{str} = \int_{t_0}^{t_f} - EA \frac{\norm{r'(l)}-1}{\norm{r'(l)}} r'(l)\cdot \delta r(l)+\int_0^l  \Big[- \mub\ddot r(\sb)\nonumber\\
&+ \mub g \,e_3 +EA \parenth{\frac{\norm{r'(\sb)}-1}{\norm{r'(\sb)}} r'(\sb)}'\Big]\cdot \delta r(\sb)\;  d\sb\, dt.\label{eqn:delG1}
\end{align}
(See \cite{KuhSeiAMC95} for details.)

Next, we found the variation of $\mathfrak{G}_{rb}$. It can be written as
\begin{align}
\delta\mathfrak{G}_{rb} = \int_{t_0}^{t_f} & \bracket{M\dot r(l)+MR\hat\Omega\rho_c}\cdot\delta \dot r(l) +Mge_3\cdot \delta r(l)\nonumber\\
&+ \bracket{J\Omega + M\hat\rho_c R^T \dot r(l)}\cdot\delta\Omega\nonumber\\
&+ M\dot r(l)\cdot \delta R\hat\Omega\rho_c + Mge_3\cdot \delta R\rho_c\;dt.\label{eqn:Grb}
\end{align}
The variation of a rotation matrix can be written as
\begin{align*}
\delta R = \frac{d}{d\epsilon}\bigg|_{\epsilon=0} R^\epsilon = \frac{d}{d\epsilon}\bigg|_{\epsilon=0} R \exp \epsilon\hat\eta = R \hat\eta 
\end{align*}
for $\eta\in\Re^3$~\cite{LeeLeoPICCA05}. The corresponding variation of the angular velocity is obtained from the kinematics equation \refeqn{Rdot}:
\begin{align*}
\delta\hat\Omega = \frac{d}{d\epsilon}\bigg|_{\epsilon=0} (R^\epsilon)^T \dot R^\epsilon = (\dot\eta + \Omega\times \eta)^\wedge.
\end{align*}
Substituting these into \refeqn{Grb} and applying the integration by parts, we obtain
\begin{align}
\delta& \mathfrak{G}_{rb} = \int_{t_0}^{t_f} -\bracket{M\ddot r(l)-MR\hat\rho_c\dot\Omega+MR\hat\Omega^2\rho_c-Mge_3}\cdot \delta r(l)\nonumber\\ 
& + \bracket{-J\dot\Omega - M\hat\rho_c R^T \ddot r(l)
+M\hat\rho_c\hat\Omega R^T \dot r(l)}\cdot\dot\eta\nonumber\\
& +\bracket{-M{\hat\rho_c\hat\Omega}R^T\dot r(l) + Mg \hat\rho_c R^T e_3
-\hat\Omega J\Omega }\cdot \eta
\;dt,\label{eqn:delG2}
\end{align}
where we repeatedly use the property: $y\cdot \hat x z = \hat z y\cdot x$ for any $x,y,z\in\Re^3$.

From \refeqn{delG1} and \refeqn{delG2}, the variation of the action integral is given by $\delta \mathfrak{G}=\delta\mathfrak{G}_{str}+\delta\mathfrak{G}_{rb}$, and it is equal to zero for any variation according to Hamilton's principle. This yields the following Euler-Lagrange equations:
\begin{gather}
\mub \ddot r(\sb,t)  - \mub g \,e_3 -EA \deriv{}{\sb}\parenth{\frac{\norm{r'(\sb,t)}-1}{\norm{r'(\sb,t)}} r'(\sb,t)}=0,\\
\begin{aligned}
M\Bigl( \ddot r(l,t) - & R\hat\rho_c  \dot\Omega +R\hat\Omega^2\rho_c-g e_3\Bigr) \\ &+ EA \frac{\norm{r'(l,t)}-1}{\norm{r'(l,t)}} r'(l,t) =0,\end{aligned}\\
J\dot\Omega+\hat\Omega J\Omega+m\hat\rho_c R^T\ddot r(l,t) -mg\hat\rho_cR^Te_3=0.
\end{gather}

\paragraph*{Conserved quantities} The total energy, given by $E=T_{str}+V_{str}+T_{rb}+V_{rb}$, is preserved. As the Lagrangian is invariant under the rotation about the gravity direction, the total angular momentum about the gravity direction is conserved. It is given by
$\pi_3 = \{ \int_{0}^l \mub \hat r(\sb) \dot r(\sb) \,d\sb + M \hat r(l) (\dot r(l)+R\hat\Omega\rho_c) -M\hat{\dot {r}} (l) R\rho_c + RJ\Omega\} \cdot e_3.$

\section{Lie Group Variational Integrator}\label{sec:LGVI}

The continuous-time Euler-Lagrange equations developed in the previous section provide an analytical model for a 3D elastic string pendulum. However, the popular finite difference approximations or finite element approximations of those equations using a general purpose numerical integrator may not preserve the geometric properties of the system accurately~\cite{HaiLub00}.

Variational integrators provide a systematic method of developing geometric numerical integrators for Lagrangian/Hamiltonian systems~\cite{MarWesAN01}. As it is derived from a discrete analogue of Hamilton's principle, it preserves symplecticity and the momentum map, and it exhibits good total energy behavior. Lie group methods conserve the structure of a Lie group configuration manifold as it updates a group element using the group operation~\cite{IseMunAN00}.

These two methods have been unified to obtain a Lie group variational integrator for Lagrangian/Hamiltonian systems evolving on a Lie group~\cite{Lee08}. This preserves symplecticity and group structure of those systems concurrently. It has been shown that this property is critical for accurate and efficient simulations of rigid body dynamics~\cite{LeeLeoCMDA07}.

In this section, we develop a Lie group variational integrator for a 3D elastic string pendulum. We first construct a finite element model, and derive an expression for a discrete Lagrangian, which is  substituted into the discrete-time Euler-Lagrange equations on a Lie group.

\subsection{Finite Element Model}

We discretize the string by $N$ one-dimensional line elements. Thus, the unstretched length of 
each element is $u = \frac{l}{N}$. A natural coordinate $\zeta\in[0,1]$ in the $a$-th element is defined by $\zeta = \frac{1}{u} (\sb - u (a-1))$. Let $S_0,S_1$ be shape functions given by $S_0(\zeta)=1-\zeta$, and $S_1(\zeta)=\zeta$. These shape functions are also referred to as  \textit{tent functions}. The position vectors for the end nodes of the $a$-th element are given by $r_{k,a},r_{k,a+1}$ when $t=kh$ for a fixed time step $h$.

Using this finite element model, the position vector $r(\sb,t)$ of a material point in the $a$-th element is approximated as follows:
\begin{align}
r(\sb,t) = S_0(\zeta) r_{k,a} + S_1(\zeta) r_{k,a+1} \equiv r_{k,a}(\zeta).\label{eqn:rFEM}
\end{align}
Note that $r_{k,a}(0)=r_{k,a}$ and $r_{k,a}(1)=r_{k,a+1}$. The partial derivative with respect to $\sb$ is given by
\begin{align}
r'(\sb,t)= \deriv{r(\sb,t)}{\zeta}\deriv{\zeta}{\sb} = 
\frac{1}{u}(r_{k,a+1}-r_{k,a})\equiv r'_{k,a}.\label{eqn:rpFEM}
\end{align}
The partial derivative with respect to $t$ is approximated by
\begin{align}
\dot r (\sb,t) = \frac{1}{h} ( S_0(\zeta) \Delta r_{k,a} + S_1(\zeta)\Delta r_{k,a+1}) \equiv v_{k,a}(\zeta),\label{eqn:vFEM}
\end{align}
where the Delta-operator represents a change for a time step, i.e. $\Delta r_{k,a}= r_{k+1,a}-r_{k,a}$. 

\subsection{Discrete-Lagrangian}

Using these finite element model, a configuration of the discretized 3D elastic pendulum at $t=kh$ is described by $g_k=(r_{k,1},\ldots,r_{k,N+1},R_k)$, and the corresponding configuration manifold is $\G=(\Re^3)^{N+1}\times \SO$. 

We define a discrete-time kinematics equation as follows. Define $f_k=(\Delta r_{k,1},\ldots,\Delta r_{k,N+1},F_k)\in\G$ for $\Delta r_{k,a}\in\Re^3$ and $F_k\in\SO$ such that $g_{k+1}=g_k f_k$ and $\G$ acts on itself by the diagonal action:
\begin{align}
    (&r_{k+1,1},\ldots,r_{k+1,N+1},R_{k+1})\nonumber\\
    & =(r_{k,1}+\Delta r_{k,1},\ldots,r_{k,N+1}+\Delta r_{k,N+1},R_k F_k).
\end{align}
Therefore, $f_k$ represents the relative update between two integration steps. This ensures that the structure of the Lie group configuration manifold is numerically preserved since $g_{k}$ is updated by $f_k$ using the right Lie group action of $\G$ on itself.

A discrete Lagrangian $L_d(g_k,f_k):\G\times\G\rightarrow\Re$ is an approximation of the Jacobi solution of the Hamilton--Jacobi equation, which is given by the integral of the Lagrangian along the exact solution of the Euler-Lagrange equations over a single time step:
\begin{align*}
    L_d(g_k,f_k)\approx \int_0^h L(\tilde g(t),{\tilde g}^{-1}(t)\dot{\tilde g} (t))\,dt,
\end{align*}
where $\tilde g(t):[0,h]\rightarrow \G$ satisfies Euler-Lagrange equations with boundary conditions $\tilde{g}(0)=g_k$, $\tilde{g}(h)=g_kf_k$. The resulting discrete-time Lagrangian system, referred to as a variational integrator, approximates the Euler-Lagrange equations to the same order of accuracy as the discrete Lagrangian approximates the Jacobi solution.

Substituting \refeqn{rFEM}-\refeqn{vFEM} into the continuous-time Lagrangian given by \refeqn{L}, the contribution of the $a$-th element to the discrete Lagrangian is chosen as follows.
\begin{align*}
L_{d_{k,a}} & = \int_{0}^1 \frac{1}{h} \mub \norm{v_{k,a}(\zeta)}^2\;u d\zeta\nonumber\\
&- \frac{h}{2} \int_0^1 \frac{1}{2} EA (\norm{r'_{k,a}}-1)^2 -\mub g\, r_{k,a}(\zeta)\cdot e_3\;u d\zeta \\ 
&- \frac{h}{2} \int_0^1 \frac{1}{2} EA (\norm{r'_{k+1,a}}-1)^2 -\mub g\, r_{k+1,a}(\zeta)\cdot e_3\;u d\zeta.
\end{align*}
This is given by
\begin{align}
L_{d_{k,a}}& = \frac{1}{6h}m\Delta r_{k,a}\cdot \Delta r_{k,a} 
+ \frac{1}{6h}m \Delta r_{k,a}\cdot \Delta r_{k,a+1} \nonumber\\
& + \frac{1}{6h}m \Delta r_{k,a+1}\cdot \Delta r_{k,a+1} \nonumber\\
&+ \frac{h}{4}m g (2r_{k,a}+2r_{k,a+1}+\Delta r_{k,a}+\Delta r_{k,a+1})\cdot e_3 \nonumber\\
&- \frac{1}{4}h \kappa ( \norm{r_{k,a+1}-r_{k,a}}-u)^2\nonumber\\
&- \frac{1}{4}h \norm{r_{k,a+1}+\Delta r_{k,a+1}-r_{k,a}-\Delta r_{k,a}}-u)^2),\label{eqn:Ldka}
\end{align}
where $m = \mub u$, $\kappa = \frac{EA}{u}$. So, the contribution of the string to the discrete Lagrangian is $L_{d_{k,str}}=\sum_{a=1}^{N} L_{d_{k,a}}$. The contribution of the rigid body to the discrete Lagrangian is chosen as follows.
\begin{align}
L_{d_{k,rb}}& =  \frac{1}{2h} M \Delta r_{k,N+1}\cdot \Delta r_{k,N+1}
+\frac{1}{h}\tr{(I-F_k)J_d}\nonumber\\
& + \frac{M}{h}\Delta r_{k,N+1}\cdot R_k(F_k-I)\rho_c\nonumber\\
& + \frac{h}{2}Mg\, (r_{k,N+1}+R_k\rho_c) \cdot e_3\nonumber\\
& + \frac{h}{2}Mg\, (r_{k,N+1}+\Delta r_{k,N+1} +R_kF_k\rho_c) \cdot e_3,\label{eqn:Ldkrb}
\end{align}
where $J_d\in\Re^{3\times 3}$ is a nonstandard inertia matrix defined by $J_d = \frac{1}{2}\tr{J}I_{3\times 3}-J$, as introduced in~\cite{LeeLeoPICCA05}. 

From \refeqn{Ldka}, \refeqn{Ldkrb}, the discrete Lagrangian of the 3D elastic string pendulum is as follows.
\begin{align}
L_{d_k}(g_k,f_k) & = L_{d_{k,str}}(g_k,f_k) + L_{d_{k,rb}}(g_k,f_k)\nonumber\\
& = \sum_{a=1}^{N} L_{d_{k,a}}(g_k,f_k) + L_{d_{k,rb}}(g_k,f_k).\label{eqn:Ldk}
\end{align}

\subsection{Discrete-time Euler-Lagrange Equations}

For a discrete Lagrangian on $\G\times\G$, the following discrete-time Euler-Lagrange equations, referred to as a Lie group variational integrator, were developed in~\cite{Lee08}.
\begin{gather}
\begin{aligned}
    \T_e^*\L_{f_{k-1}}\cdot \D_{f_{k-1}}L_{d_{k-1}}-&\Ad^*_{f_{k}^{-1}}\cdot(\T_e^*\L_{f_{k}}\cdot \D_{f_{k}}L_{d_{k}})\\
    &\quad+\T_e^*\L_{g_{k}}\cdot \D_{g_{k}} L_{d_{k}}=0,
\end{aligned}\label{eqn:DEL_G0}\\
    g_{k+1} = g_k f_k,\label{eqn:DEL_G1}
\end{gather}
where $\T\L:\T\G\rightarrow\T\G$ is the tangential map of the left translation, $\D_f$ represents the derivative with respect to $f$, and $\Ad^*:\G\times\g^*\rightarrow\g^*$ is $\mathrm{co}$-$\mathrm{Ad}$ operator~\cite{MarRat99}.

Using this result, we develop a Lie group variational integrator for a 3D elastic string pendulum. For $f=(\Delta r_1,\ldots,\Delta r_{N+1},F)\in \G$ and $p=(p_1,\ldots,p_{N+1},\pi)\in\g^*\simeq(\Re^{3})^{N+1}\times \Re^3$, the $\mathrm{co}$-$\mathrm{Ad}$ operator is given by $\Ad^*_{f^{-1}} p = (p_1,\ldots,p_{N+1},F\pi)$.

\paragraph*{Derivatives of the discrete Lagrangian}

We now obtain expressions for the derivatives of the discrete Lagrangian. The derivatives of the discrete Lagrangian of the $a$-th element, given by \refeqn{Ldka}, with respect to $\Delta r_{k,a}$ and $\Delta r_{k,a+1}$ are given by
\begin{align}
\D_{\Delta r_{k,a}} L_{d_{k,a}} & = \frac{1}{3h} m (\Delta r_{k,a} + \frac{1}{2} \Delta r_{k,a+1})
+\frac{h}{4}m g e_3\nonumber\\ &\quad + \frac{h}{2}\nabla V^e_{k+1,a},\nonumber\\
\D_{\Delta r_{k,a+1}} L_{d_{k,a}} & = \frac{1}{3h} m (\Delta r_{k,a+1} + \frac{1}{2} \Delta r_{k,a})
+\frac{h}{4}m g e_3\nonumber\\ &\quad - \frac{h}{2}\nabla V^e_{k+1,a}.\label{eqn:DdrkapLd}
\end{align}
where $\nabla V^e_{k,a} = \kappa \frac{\norm{ x}-u}{\norm { x}}  x$ for $ x=r_{k,a+1}-r_{k,a} \in \Re^3$. Then, from \refeqn{Ldk}, the derivative of the discrete Lagrangian with respect to $\Delta r_{k,a}$, for $a\in\{2,\ldots,N\}$,  is given by
\begin{align}
\D_{\Delta r_{k,a}} & L_{d_{k}} = \D_{\Delta r_{k,a}} L_{d_{k,a}} + \D_{\Delta r_{k,a}} L_{d_{k,a-1}}\nonumber\\
& = \frac{1}{6h} m (\Delta r_{k,a-1} +4\Delta r_{k,a} + \Delta r_{k,a+1})\nonumber\\
& \quad+\frac{h}{2}m g e_3 + \frac{h}{2}\nabla V^e_{k+1,a} - \frac{h}{2}\nabla V^e_{k+1,a-1}.\label{eqn:DdrLd}
\end{align}
Similarly, the derivative of the discrete Lagrangian with respect to $r_{k,a}$, for $a\in\{2,\ldots,N\}$,  is given by
\begin{align}
\D_{r_{k,a}} L_{d_{k}} & = hm g e_3 +\frac{h}{2} (\nabla V^e_{k,a}+\nabla V^e_{k+1,a})\nonumber\\
& \quad -\frac{h}{2} (\nabla V^e_{k,a-1}+\nabla V^e_{k+1,a-1}).
\end{align}

Next, we find the derivatives of the discrete Lagrangian with respect to $\Delta r_{k,N+1}$ and $r_{k,N+1}$. They are contributed by the $N$-th string element and the rigid body, and they can be obtained from \refeqn{Ldkrb} and \refeqn{DdrkapLd}  as follows.
\begin{align}
&\D_{\Delta r_{k,N+1}} L_{d_k} = \frac{1}{h} (M+\frac{m}{3}) \Delta r_{k,N+1} +\frac{1}{6h}m\Delta r_{k,N}\nonumber\\
&+\frac{M}{h}R_k(F_k-I)\rho_c
+\frac{h}{2}(M+\frac{m}{2})ge_3-\frac{h}{2}\nabla V^e_{k+1,N},\\
&\D_{r_{k,N+1}} L_{d_k} = h(M+\frac{m}{2}) g e_3 -\frac{h}{2}\nabla V^e_{k,N}-\frac{h}{2}\nabla V^e_{k+1,N}.
\end{align}

Now, we find the derivatives of the discrete Lagrangian with respect to $F_k$ and $R_k$. From \refeqn{Ldkrb}, we have
\begin{align*}
\D_{F_k} L_{d_k}\cdot \delta F_k & = \frac{1}{h}\tr{-\delta F_k J_d} + \frac{M}{h} \Delta r_{k,N+1}\cdot R_k\delta F_k\rho_c\\
&\quad + \frac{h}{2}Mg R_k\delta F_k\rho_c \cdot e_3\\
& = \frac{1}{h}\tr{-\delta F_k J_d} + A_k\cdot \delta F_k\rho_c,
\end{align*}
where $A_k = \frac{M}{h} R^T_k \Delta r_{k,N+1} + \frac{h}{2}Mg R_k^T e_3$.
The variation of $F_{k}$ can be written as $\delta F_{k}=F_{k}\hat\zeta_{k}$ for $\zeta_{k}\in\Re^3$. Therefore, this can be written as
\begin{align*}
\D_{F_{k}}L_{d_k} \cdot (F_{k}\hat\zeta_{k})
& =  (\T_I^* \L_{F_{k}}\cdot \D_{F_{k}}L_{d_k}) \cdot \zeta_{k}\\
& = \frac{1}{h}\tr{- F_{k} \hat\zeta_{k} J_{d}}
+A_k\cdot F_k \hat\zeta_k \rho_c.
\end{align*}
By repeatedly applying the following property of the trace operator, $\mbox{tr}[AB]=\mbox{tr}[BA]=\mbox{tr}[A^TB^T]$ for any $A,B\in\Re^{3\times 3}$, the first term can be written as $\mbox{tr}[-F_{k} \hat\zeta_{k} J_{d}]= \mbox{tr}[- \hat\zeta_{k} J_{d}F_{k}]=\mbox{tr}[\hat\zeta_{k}F_{k}^TJ_{d}]=-\frac{1}{2}\mbox{tr}[\hat\zeta_{k}(J_{d}F_{0}-F_{k}^T J_{d})]$. Using the property of the hat map, $x^T y = -\frac{1}{2}\mbox{tr}[\hat x\hat y]$ for any $x,y\in\Re^3$, this can be further written as $((J_{d}F_{k}-F_{k}^T J_{d})^\vee) \cdot \zeta_{k}$. As $y\cdot\hat x z = \hat z y\cdot x$ for any $x,y,z\in\Re^3$, the second term can be written as $F_k^T A_k\cdot \hat\zeta_k \rho_c=\hat\rho_c F_k^TA_k\cdot \zeta_k$. Using these, we obtain
\begin{align}
\T_I^* & \L_{F_{k}}\cdot \D_{F_{k}}L_{d_k}  = \frac{1}{h}(J_{d}F_{k}-F_{k}^T J_{d})^\vee+\hat\rho_c F_k^T A_k.\label{eqn:DFLd}
\end{align}
The the co-$\Ad$ operator yields
\begin{align}
\Ad^*_{F_k^T}\cdot(\T_I^* & \L_{F_{k}}\cdot \D_{F_{k}}L_{d_k})  = \frac{1}{h}(F_kJ_{d}- J_{d}F_{k}^T)^\vee+\widehat {F_k\rho_c} A_k.\label{eqn:AdDFLd}
\end{align}

Similarly, we can derive the derivative of the discrete Lagrangian with respect to $R_k$ as follows.
\begin{align}
\T_I^* \L_{R_{k}}\cdot \D_{R_{k}}L_{d_k}&=\frac{M}{h} ((F_k-I)\rho_c)^\wedge R_k^T \Delta r_{k,N+1}\nonumber\\&\quad+\frac{h}{2}Mg\hat\rho_c R_k^T e_3+\frac{h}{2}Mg \widehat{F_k\rho_c} R_k^T e_3.\label{eqn:DRLd}
\end{align}

\paragraph*{Discrete-time Euler-Lagrange Equations} Substituting \refeqn{DdrLd}-\refeqn{DRLd} into \refeqn{DEL_G0}-\refeqn{DEL_G1}, we obtain discrete-time Euler-Lagrange equations for a 3D elastic string pendulum as follows.
\begin{gather}
\begin{aligned}
\frac{1}{6h} m & (\Delta^2 r_{k,a-1}+4\Delta^2 r_{k,a}+\Delta^2 r_{k,a+1})\\
&- hm g e_3
+h\nabla V^e_{k,a-1}
-h\nabla V^e_{k,a}=0,
\end{aligned}\label{eqn:DEL0}\\
\begin{aligned}
&\frac{1}{h} (M+\frac{m}{3}) \Delta^2 r_{k,N+1}+\frac{1}{6h} m \Delta^2 r_{k,N}+h\nabla V^e_{k,N}
\\
&+\frac{1}{h}M(R_kF_k-2R_k+R_{k-1})\rho_c
-h(M+\frac{m}{2})ge_3
=0,
\end{aligned}\label{eqn:DEL1}\\
\begin{aligned}
\frac{1}{h}(F_kJ_d&-J_dF_k^T-J_dF_{k-1}+F_{k-1}^TJ_d)^\vee\\ 
&+\frac{M}{h}\hat\rho_c R_{k}^T \Delta^2 r_{k,N+1}-hMg\hat\rho_c R_k^T e_3=0,
\end{aligned}\label{eqn:DEL2}\\
r_{k+1,a}= r_{k,a}+\Delta r_{k,a},\label{eqn:DEL3}\\
R_{k+1}=R_{k}F_k.\label{eqn:DEL4}
\end{gather}
where $\Delta^2 r_{k,a} = \Delta r_{k,a}-\Delta r_{k-1,a} = r_{k+1,a} - 2r_{k,a}+r_{k-1,a}$, $u=\frac{l}{N}$, $m=\mub u$, $\kappa=\frac{EA}{u}$, and $\nabla V^e_{k,a} = \kappa \frac{\norm{ x}-u}{\norm { x}}  x$ for $x=r_{k,a+1}-r_{k,a}$. Equation \refeqn{DEL0} is satisfied for $a\in\{2,\ldots,N\}$, and \refeqn{DEL3} is satisfied for $a\in\{2,\ldots,N+1\}$. For any $k$, the vector $r_{k,1}=0$ since the pivot is fixed.

For given $(g_{k-1},f_{k-1})$, $g_k$ is explicitly computed by \refeqn{DEL3} and \refeqn{DEL4}. The update $f_k$ is computed by a fixed point iteration for $F_k$: we select an initial guess of $F_k$; $\Delta r_{k,a}$ is obtained by solving \refeqn{DEL0} and \refeqn{DEL1}, which requires the inversion of a fixed $3N\times 3N$ matrix; a new $F_k$ is computed by solving the implicit equation \refeqn{DEL2}; these are repeated until $F_k$ converges. When solving the implicit equation \refeqn{DEL2}, we first express $F_k$ on $\Re^3$ using the Cayley transform, and apply Newton's iteration~(See Section 3.3.8 in \cite{Lee08}). These yields a Lagrangian flow map $(g_{k-1},f_{k-1})\mapsto(g_k,f_k)$, and they are repeated.

\section{Numerical Example}\label{sec:NE}

We now demonstrate the computational properties of the Lie group variational integrator developed in the previous section by considering a numerical example. The material properties of the string are chosen to represent a rubber string as follows~\cite{KuhSeiAMC95}.
\begin{gather*}
\mub=0.025\,\mathrm{kg/m},\quad
l=1\,\mathrm{m},\quad EA=40\,\mathrm{N}.
\end{gather*}
The rigid body is chosen as an elliptic cylinder with a semimajor axis $0.06\,\mathrm{m}$, a semiminor axis $0.04\,\mathrm{m}$, and a height $0.1\,\mathrm{m}$. Its properties are as follows.
\begin{gather*}
M=0.1\,\mathrm{kg},\quad \rho_c= [0.04,0.01,0.05]\,\mathrm{m},\\
\quad J = \begin{bmatrix}
    0.38&   -0.04&   -0.20\\
   -0.04&    0.58&   -0.05\\
   -0.20&   -0.05&    0.30\end{bmatrix}\,\mathrm{kg\,m^2}.
\end{gather*}
Initially, the string is aligned to the horizontal $e_1$ axis at rest, and the rigid body has an initial velocity $[0,0.2,-0.5]\,\mathrm{m/s}$. We use $N=20$ elements. Simulation time is $T=5$ seconds, and time step is $h=0.0001$ second.

\begin{figure}
\centerline{
    \subfigure[$t\in{[0,1.25]}$]{
        \includegraphics[width=0.5\columnwidth]{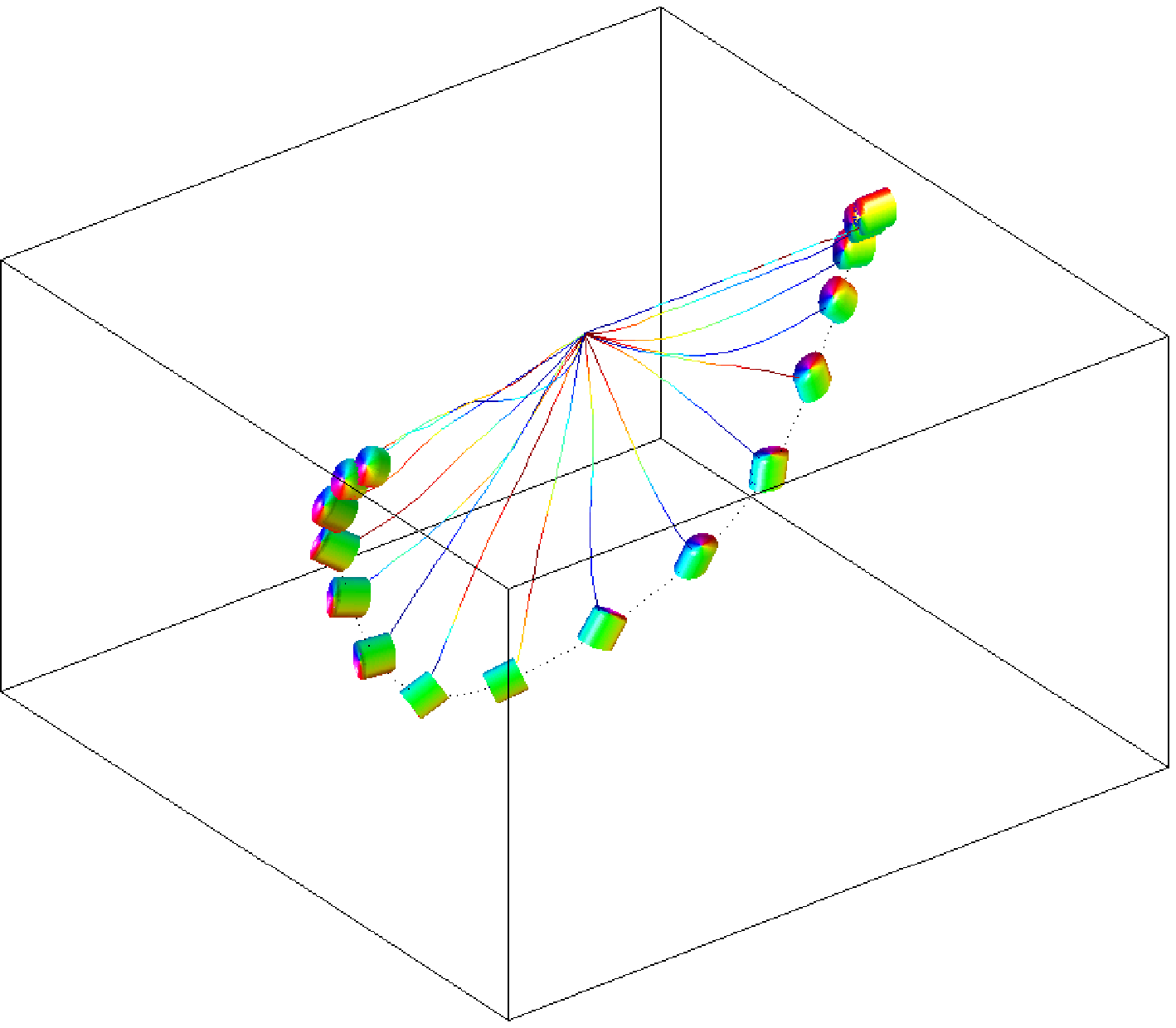}}
    \subfigure[$t\in{[1.25,2.5]}$]{
        \includegraphics[width=0.5\columnwidth]{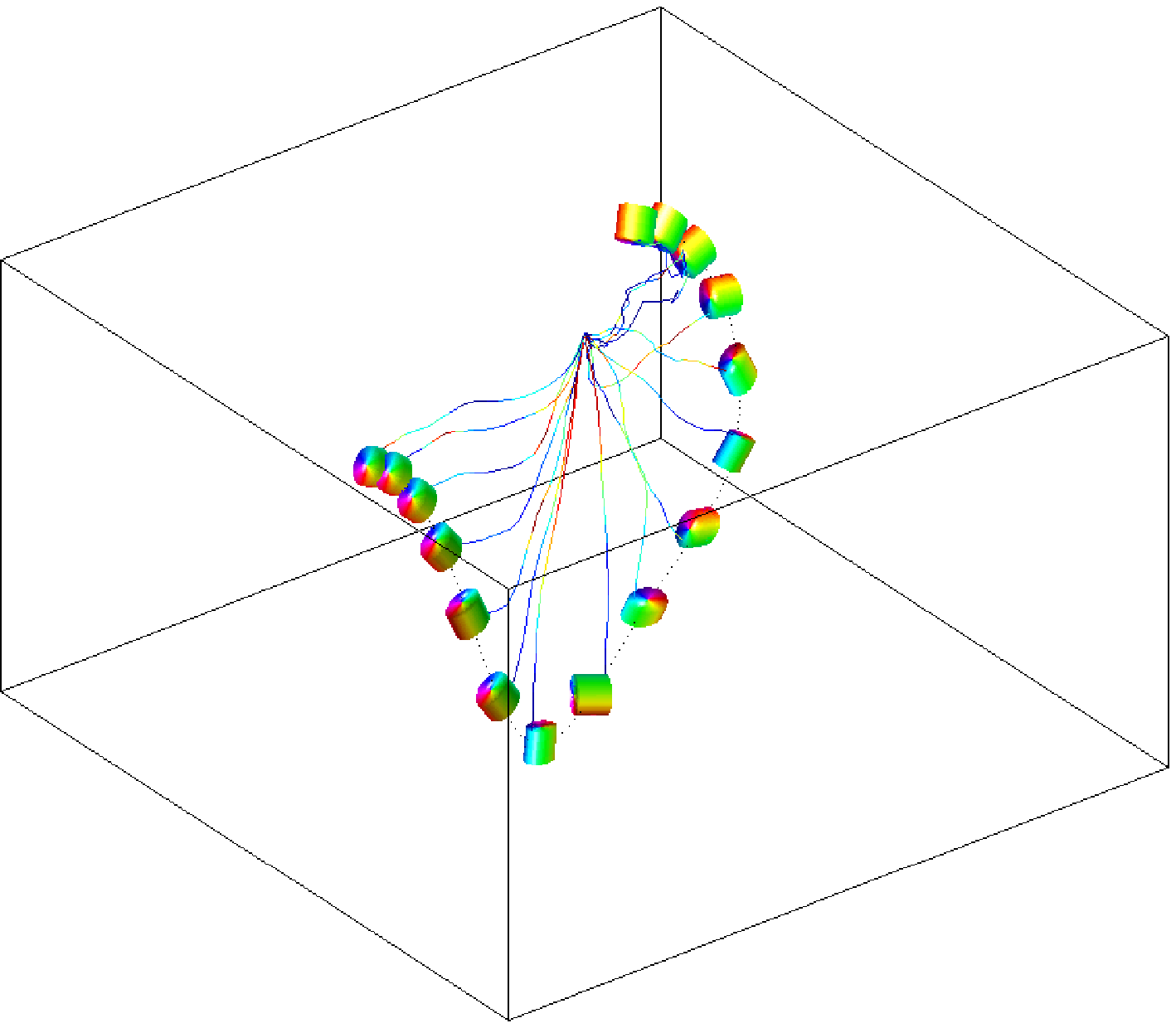}}
}
\centerline{
    \subfigure[$t\in{[2.5,3.75]}$]{
        \includegraphics[width=0.5\columnwidth]{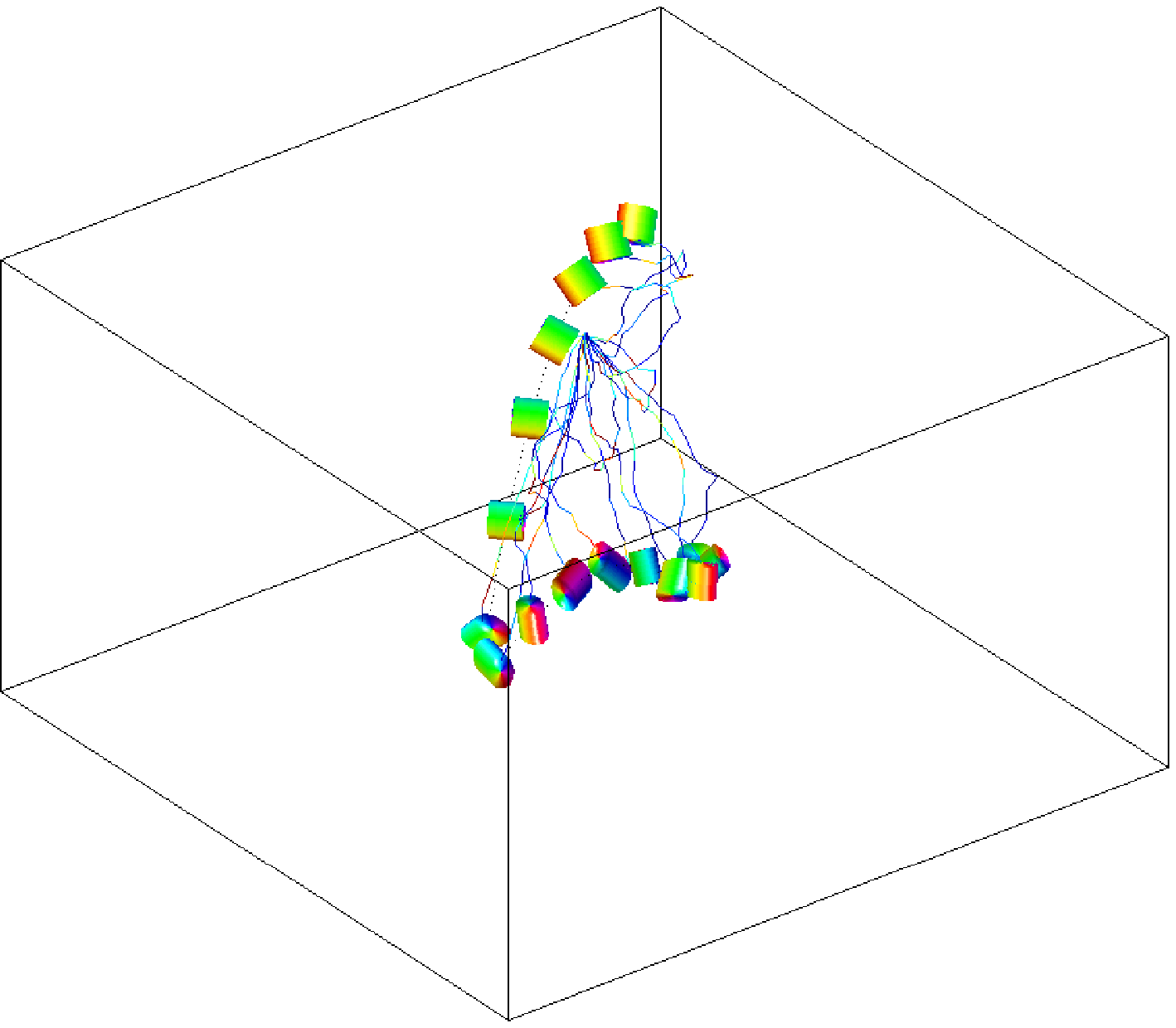}}
    \subfigure[$t\in{[3.75,5]}$]{
        \includegraphics[width=0.5\columnwidth]{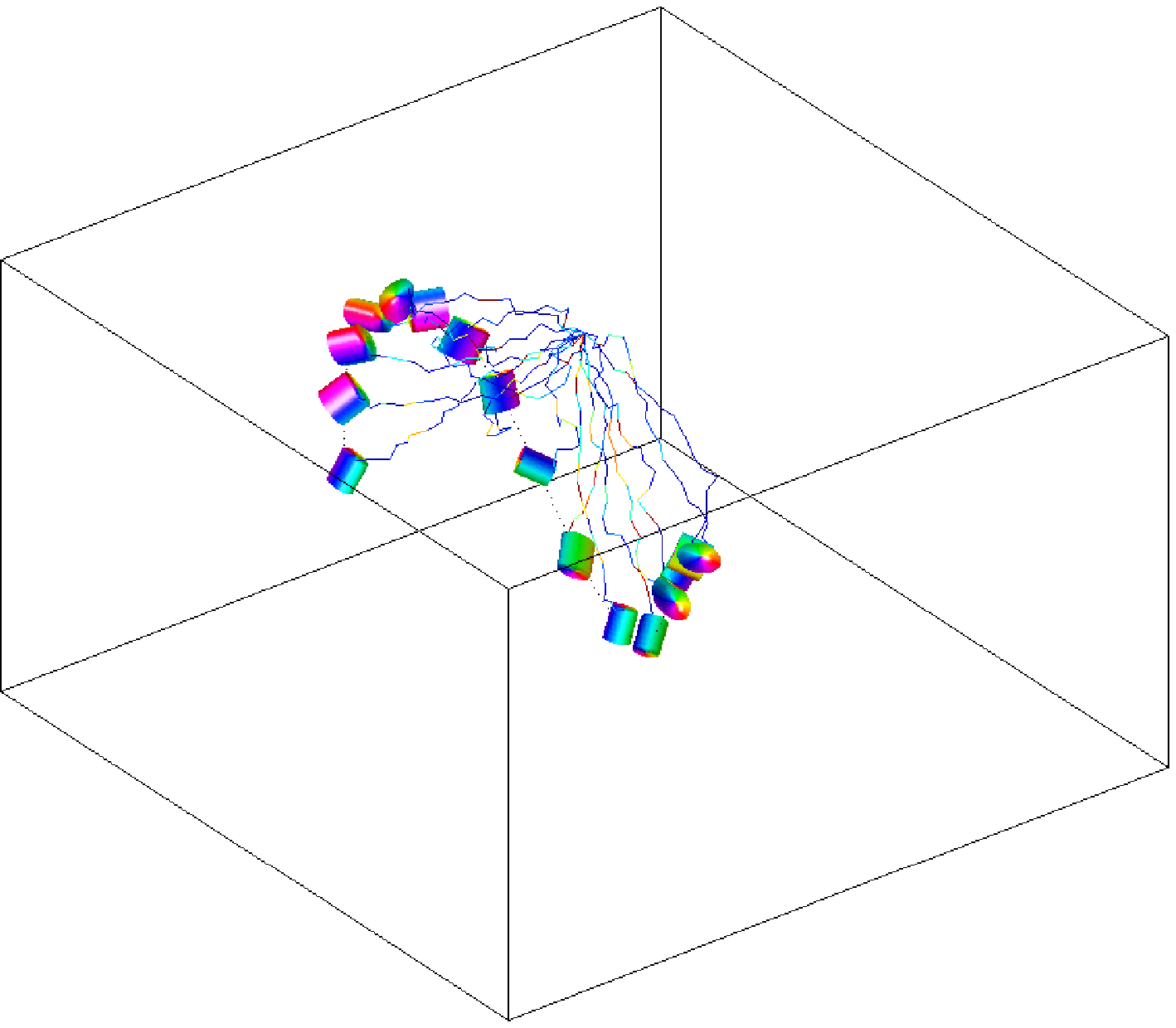}}
}
\caption{Snapshots of a 3D elastic string pendulum maneuver. Strain energy distribution is illustrated by color shading (An animation is available at http://my.fit.edu/\~{}taeyoung)}\label{fig:SM}
\end{figure}

\reffig{SM} illustrates the resulting maneuver of the 3D elastic string pendulum. As the point where the string is attached to the rigid body is displaced from the center of mass of the rigid body, the rigid body dynamics are directly coupled to the elastic string dynamics, which yields the illustrated complex maneuver.

\reffig{NS} shows the corresponding energy transfer, total energy, total angular momentum deviation, orthogonality errors of rotation matrices, velocities of the rigid body, and the stretched length of the string. As shown in \reffig{E}, the computed total energy of the Lie group variational integrator oscillates near the initial value, but there is no increasing or decreasing drift for long time periods. This is due to the fact that the numerical solutions of symplectic numerical integrators are exponentially close to the exact solution of a perturbed Hamiltonian~\cite{HaiANM94}. The value of the perturbed Hamiltonian is preserved in the discrete-time flow. The Lie group variational integrator preserves the momentum map as in \reffig{pi3}, and it also preserves the orthogonal structure of rotation matrices accurately. The orthogonality errors, measured by $\|I-R^T R\|$, are less than $2\times 10^{-13}$  in \reffig{errR}. 

These show that the Lie group variational integrator preserves the geometric characteristic of the 3D elastic string pendulum accurately even for the presented complex maneuver that has nontrivial energy transfer between different dynamic modes.

\begin{figure}
\centerline{
    \subfigure[Energy transfer (kinetic energy of the rigid body: red, kinetic energy of the string: black, gravitational potential: green, elastic potential: blue)]{
        \includegraphics[width=0.492\columnwidth]{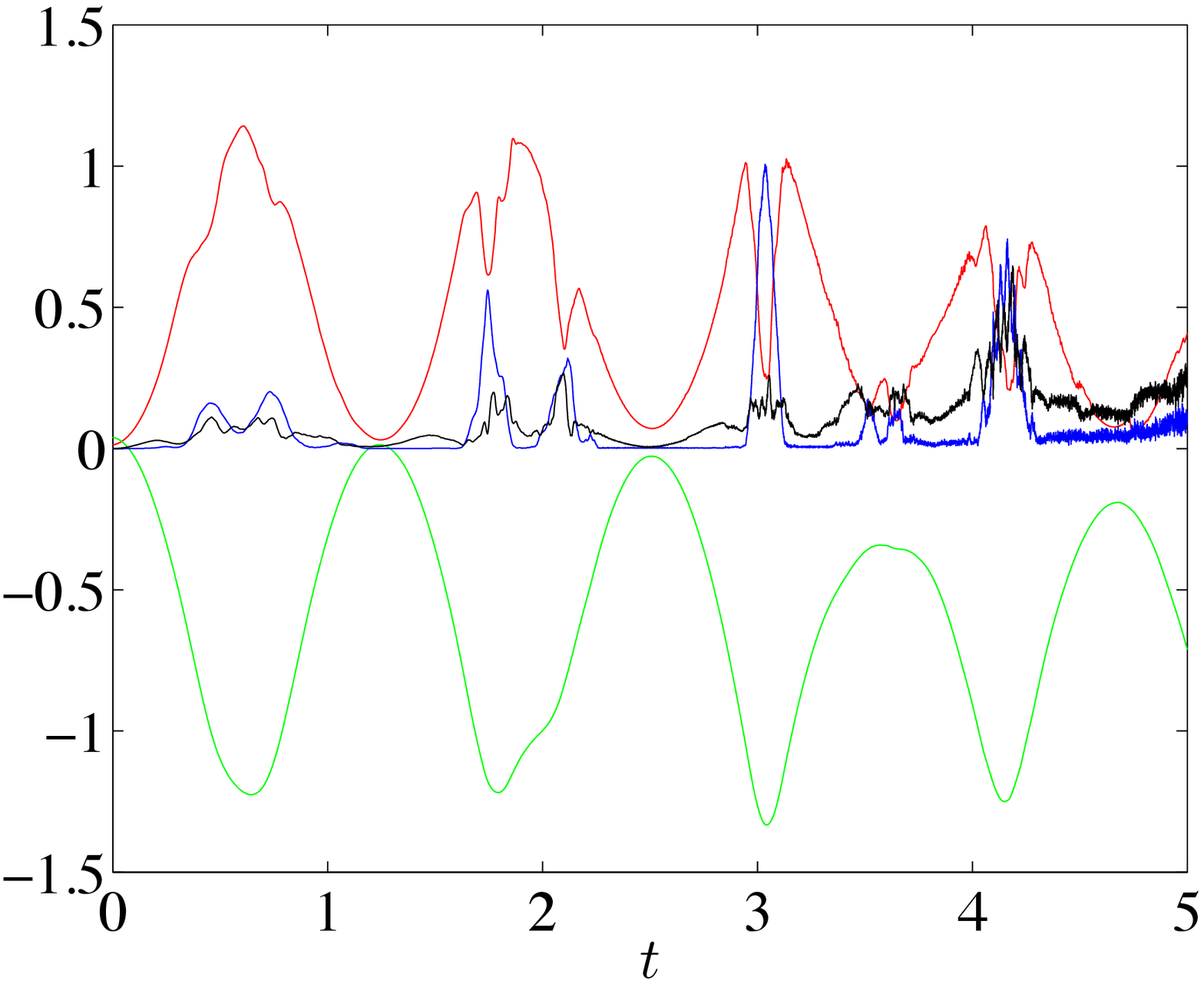}}
    \subfigure[Total energy]{
        \includegraphics[width=0.50\columnwidth]{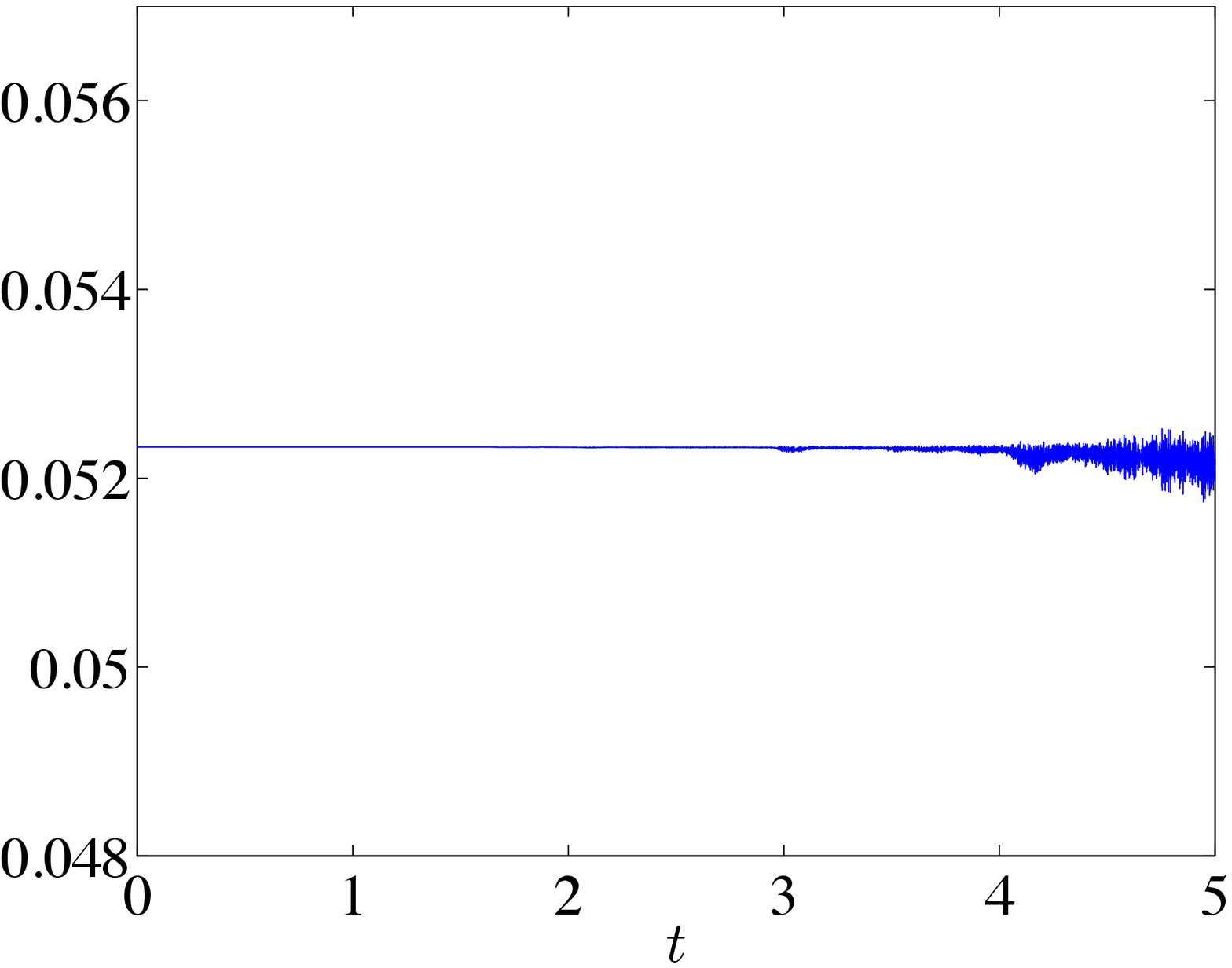}\label{fig:E}}
}
\centerline{
    \subfigure[Deviation of the total angular momentum about the gravity direction]{
        \includegraphics[width=0.50\columnwidth]{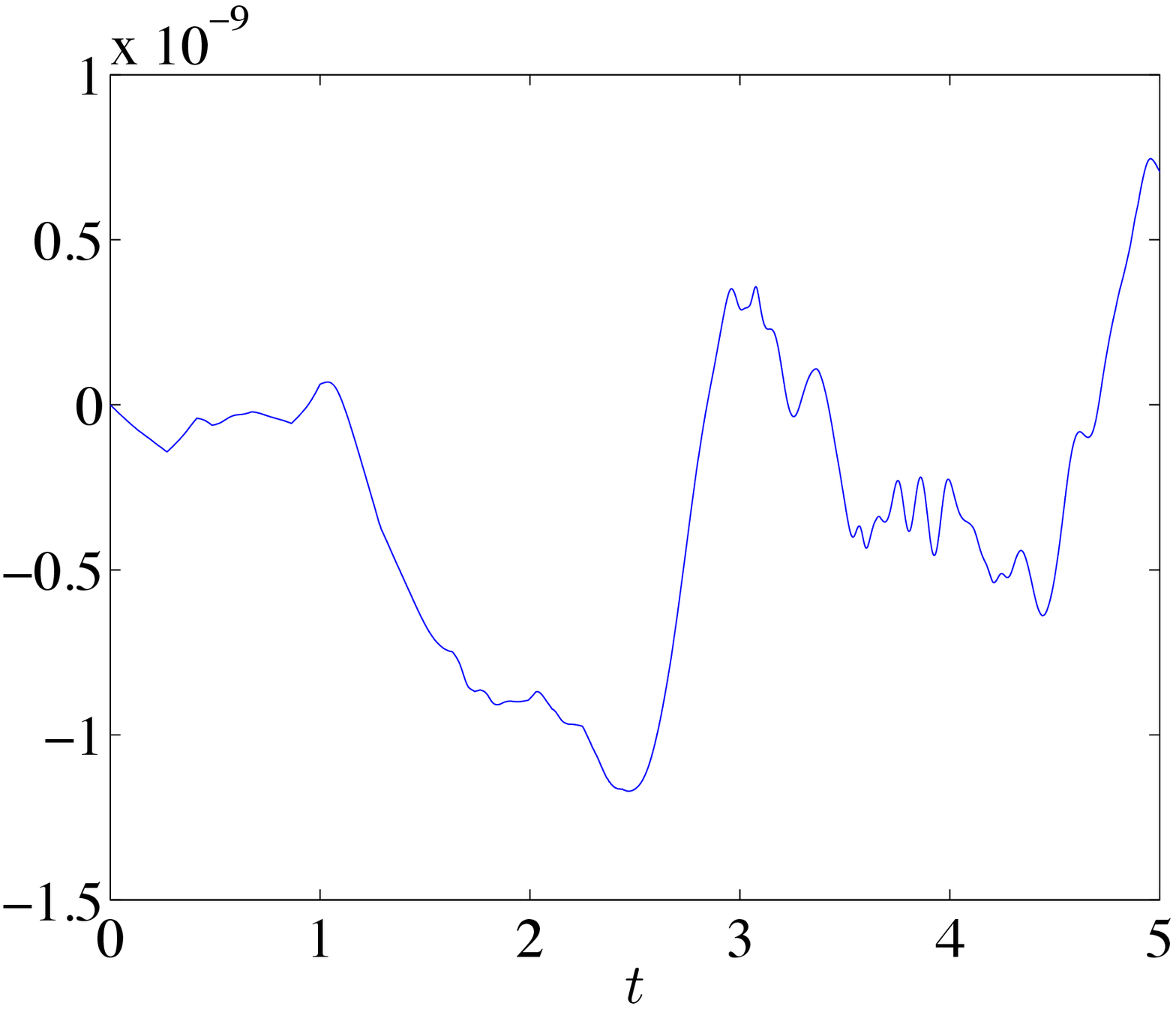}\label{fig:pi3}}
        \hspace*{0.005\textwidth}
    \subfigure[Orthogonality error of rotation matrices $\|I-R^TR\|$]{
        \includegraphics[width=0.48\columnwidth]{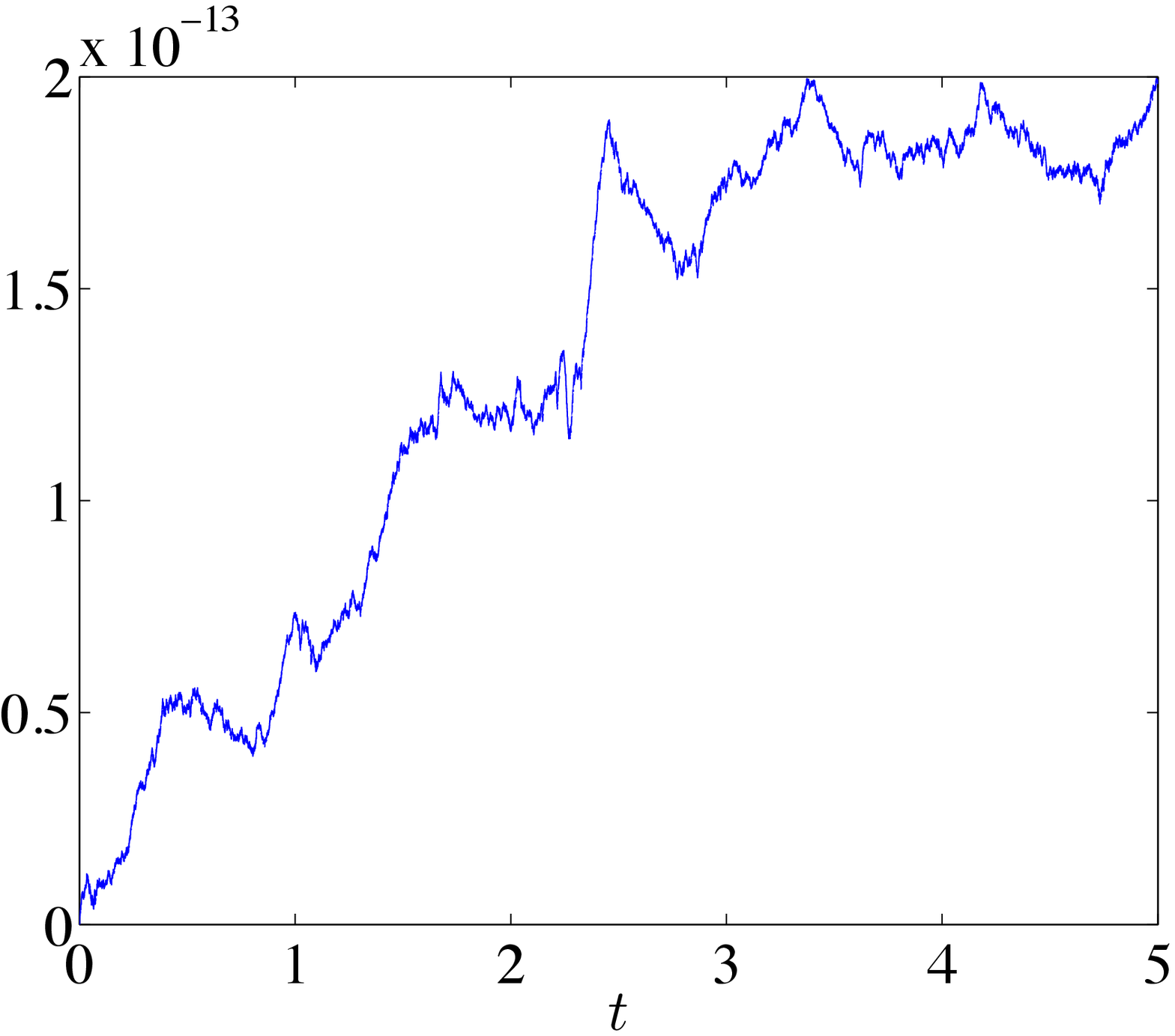}\label{fig:errR}}
}
\centerline{
    \subfigure[Velocity / angular velocity of the rigid body (second components)]{
        \includegraphics[width=0.50\columnwidth]{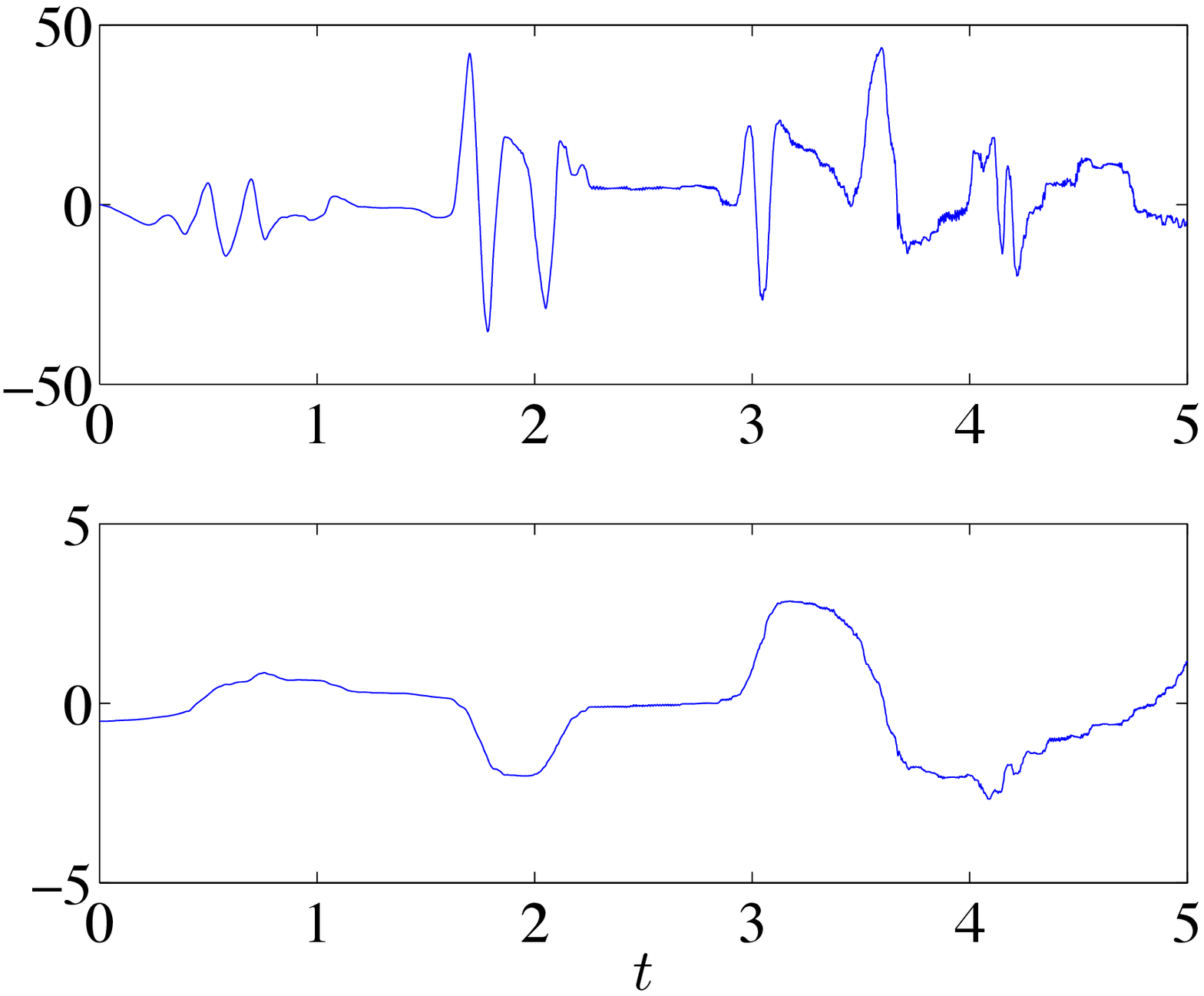}}
    \subfigure[Stretched length of the string]{
        \includegraphics[width=0.50\columnwidth]{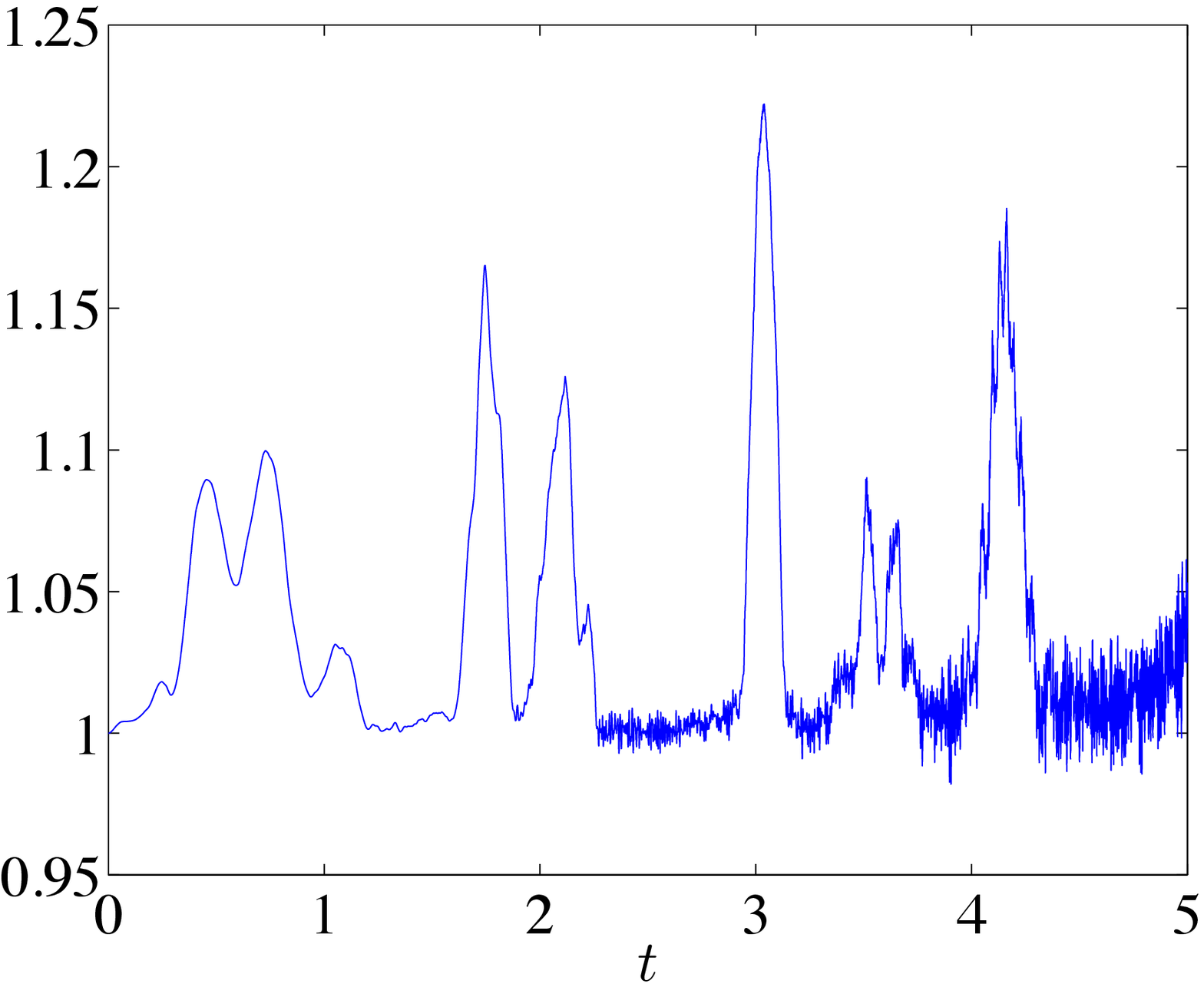}}
}
\caption{Numerical simulation of a 3D elastic string pendulum}\label{fig:NS}
\end{figure}

\section{Conclusions}

We have developed continuous-time equations of motion and a geometric numerical integrator, referred to as a Lie group variational integrator, for a 3D elastic string pendulum. The continuous-time equations of motion provide an analytical model that is defined globally on the Lie group configuration manifold, and the Lie group variational integrator preserves the geometric features of the system, thereby yielding a reliable numerical simulation tool for complex maneuvers over a long time period.

\enlargethispage{-3.5in}

These can be extended to include the effects of control inputs by using the discrete Lagrange-d'Alembert principle \cite{KanMarIJNME00}, which modifies the discrete Hamilton's principle by taking into account the virtual work of the external control inputs. When applied to an optimal control problem, this allows us to find optimal maneuvers accurately and efficiently, as there is no artificial numerical dissipation induced by the computational method. Furthermore, optimal large-angle rotational maneuvers can be easily obtained without singularities and complexity associated with local parameterizations, since the configuration is represented globally on the Lie group~\cite{LeeLeoCISSIDRWB09}.

\bibliography{CDC09}

\begin{thebibliography}{10}
\providecommand{\url}[1]{#1}
\csname url@rmstyle\endcsname
\providecommand{\newblock}{\relax}
\providecommand{\bibinfo}[2]{#2}
\providecommand\BIBentrySTDinterwordspacing{\spaceskip=0pt\relax}
\providecommand\BIBentryALTinterwordstretchfactor{4}
\providecommand\BIBentryALTinterwordspacing{\spaceskip=\fontdimen2\font plus
\BIBentryALTinterwordstretchfactor\fontdimen3\font minus
  \fontdimen4\font\relax}
\providecommand\BIBforeignlanguage[2]{{%
\expandafter\ifx\csname l@#1\endcsname\relax
\typeout{** WARNING: IEEEtran.bst: No hyphenation pattern has been}%
\typeout{** loaded for the language `#1'. Using the pattern for}%
\typeout{** the default language instead.}%
\else
\language=\csname l@#1\endcsname
\fi
#2}}

\bibitem{McLRocPI92}
T.~McLain and S.~Rock, ``Experimental measurement of rov tether tension,'' in
  \emph{Proceedings of Intervention/ROV 92}, 1992.

\bibitem{ChaOE84}
D.~Chapman, ``Towed cable behaviour during ship turning manoeuvers,''
  \emph{Ocean Engineering}, vol.~11, no.~4, pp. 327--361, 1984.

\bibitem{DriLueAOR00}
R.~Driscoll, R.~Lueck, and M.~Nahon, ``Development and validation of a
  lumped-mass dynamics model of a deep sea {R}{O}{V} system,'' \emph{Applied
  Ocean Research}, vol.~22, no.~3, pp. 169--182, 2000.

\bibitem{WalPolMC60}
T.~Walton and H.~Polacheck, ``Calculation of transient motion of submerged
  cables,'' \emph{Mathematics of Computation}, vol.~14, no.~69, pp. 27--46,
  1960.

\bibitem{KohZhaJEM99}
C.~Koh, Y.~Zhang, and S.~Quek, ``Low-tension cable dynamics: Numerical and
  experimental studies,'' \emph{Journal of Engineering Mechanics}, vol. 125,
  no.~3, pp. 347--354, 1999.

\bibitem{KuhSeiAMC95}
A.~Kuhn, W.~Seiner, J.~Zemann, D.~Dinevski, and H.~Troger, ``A comparison of
  various mathematical formulations and numerical solution methods for the
  large amplitude oscillations of a string pendulum,'' \emph{Applied
  Mathematics and Computation}, vol.~67, pp. 227--264, 1995.

\bibitem{BucDriTA04}
B.~Buckham, F.~Driscoll, and M.~Nahon, ``Development of a finite element cable
  model for use in low-tension dynamics simulation,'' \emph{Journal of Applied
  Mechanics}, vol.~71, pp. 476--485, 2004.

\bibitem{BanKanJAS82}
A.~Banerjee and T.~Kane, ``Tether deployment dynamics,'' \emph{The Journal of
  the Astronautical Sciences}, pp. 347--365, 1982.

\bibitem{SchSteISDSTA00}
M.~Schagerl, A.~Steindl, and H.~Troger, ``Dynamic analysis of the deployment
  process of tethered satellite systems,'' in \emph{IUTAM-IASS Symposium on
  Deployable Structure: Theory and Application}, 2000, pp. 345--354.

\bibitem{SheSanPICDC04}
J.~Shen, A.~Sanyal, N.~Chaturvedi, D.~Bernstein, and N.~H. McClamroch,
  ``Dynamics and control of a {3D} pendulum,'' in \emph{Proceedings of IEEE
  Conference on Decision and Control}, Dec. 2004, pp. 323--328.

\bibitem{HaiLub00}
E.~Hairer, C.~Lubich, and G.~Wanner, \emph{Geometric numerical integration},
  ser. Springer Series in Computational Mechanics 31.\hskip 1em plus 0.5em
  minus 0.4em\relax Springer, 2000.

\bibitem{LeeLeoCMDA07}
T.~Lee, M.~Leok, and N.~H. McClamroch, ``Lie group variational integrators for
  the full body problem in orbital mechanics,'' \emph{Celestial Mechanics and
  Dynamical Astronomy}, vol.~98, no.~2, pp. 121--144, June 2007.

\bibitem{Lee08}
T.~Lee, ``Computational geometric mechanics and control of rigid bodies,''
  Ph.D. dissertation, University of Michigan, 2008.

\bibitem{LeLeMc2008a}
T.~Lee, M.~Leok, and N.~McClamroch, ``Optimal attitude control of a rigid body
  using geometrically exact computations on {SO(3)},'' \emph{Journal of
  Dynamical and Control Systems}, vol.~14, no.~4, pp. 465--487, 2008.

\bibitem{LeeLeoPICCA05}
T.~Lee, M.~Leok, and N.~H. McClamroch, ``A {L}ie group variational integrator
  for the attitude dynamics of a rigid body with application to the 3{D}
  pendulum,'' in \emph{Proceedings of the {IEEE} {C}onference on {C}ontrol
  {A}pplication}, 2005, pp. 962--967.

\bibitem{MarWesAN01}
J.~Marsden and M.~West, ``Discrete mechanics and variational integrators,'' in
  \emph{Acta Numerica}.\hskip 1em plus 0.5em minus 0.4em\relax Cambridge
  University Press, 2001, vol.~10, pp. 317--514.

\bibitem{IseMunAN00}
A.~Iserles, H.~Munthe-Kaas, S.~{N\o rsett}, and A.~Zanna, ``Lie-group
  methods,'' in \emph{Acta Numerica}.\hskip 1em plus 0.5em minus 0.4em\relax
  Cambridge University Press, 2000, vol.~9, pp. 215--365.

\bibitem{MarRat99}
J.~Marsden and T.~Ratiu, \emph{Introduction to Mechanics and Symmetry},
  2nd~ed., ser. Texts in Applied Mathematics.\hskip 1em plus 0.5em minus
  0.4em\relax Springer-Verlag, 1999, vol.~17.

\bibitem{HaiANM94}
E.~Hairer, ``Backward analysis of numerical integrators and symplectic
  methods,'' \emph{Annals of Numerical Mathematics}, vol.~1, no. 1-4, pp.
  107--132, 1994, scientific computation and differential equations (Auckland,
  1993).

\bibitem{KanMarIJNME00}
C.~Kane, J.~Marsden, M.~Ortiz, and M.~West, ``Variational integrators and the
  {N}ewmark algorithm for conservative and dissipative mechanical systems,''
  \emph{International Journal for Numerical Methods in Engineering}, vol.~49,
  no.~10, pp. 1295--1325, 2000.

\bibitem{LeeLeoCISSIDRWB09}
\BIBentryALTinterwordspacing
T.~Lee, M.~Leok, and N.~H. McClamroch, ``Computational geometric optimal
  control of rigid bodies,'' \emph{Communications in Information and Systems,
  special issue dedicated to R. W. Brockett}, 2009, accepted. [Online].
  Available: \url{http://arxiv.org/abs/0805.0639}
\BIBentrySTDinterwordspacing

\end{thebibliography}
\bibliographystyle{IEEEtran}

\end{document}